\documentclass[letterpaper, 10pt, conference]{ieeeconf}

\IEEEoverridecommandlockouts
\usepackage{amsmath}
\usepackage{mynotation}
\usepackage{pict2e,picture}
\usepackage[abbr=false]{clevethm}
\usepackage{algorithm}
\usepackage[noend]{algpseudocode}
\usepackage{hyperref}

\definecolor{linkblue}{RGB}{11,104,149}
\hypersetup{
    colorlinks=true,
    linkcolor=linkblue,
    urlcolor=linkblue,
}

\usepackage{xcolor}
\usepackage{pgfplots}
\usepackage{pgfplotstable}
\usepackage{float}
\pgfplotsset{compat=1.7}

\newif\ifArxiv
\Arxivfalse

\newif\ifMergefig
\Mergefigtrue

\newcommand{\algsize}{\footnotesize}

\pgfplotsset{
SmallBarPlot/.style={
    font=\small,
    ybar,
    width=\linewidth,
    ymin=0,
    xtick=data,
    bar width=0.4em,
},
}

\let\labelindent\relax  
\usepackage{enumitem}
\newlist{inlinelist}{enumerate*}{1}
      \setlist*[inlinelist,1]{%
        label=(\arabic*),
      }
      \newlist{inlinelist*}{enumerate*}{1}
      \setlist[inlinelist*,1]{%
        label={(\arabic*)},
        itemjoin={;\ },itemjoin* = {;\ and },
}

\usepackage{acro}
\DeclareAcronym{nmpc}{
  short=NMPC,
  long=nonlinear model predictive control
}

\usepackage[style=ieee]{biblatex}

\newcommand{\reals}{\ensuremath{\Re}}
\newcommand{\sK}{\ensuremath{\mathcal{K}}}
\newcommand{\sJ}{\ensuremath{\mathcal{J}}}

\newcommand{\hessxx}{\ensuremath{\nabla^2_{\!xx}\!}}

\makeatletter
\DeclareRobustCommand{\bigtimes}{%
  \mathop{\vphantom{\sum}\mathpalette\@bigtimes\relax}\slimits@
}
\newcommand{\@bigtimes}[2]{\vcenter{\hbox{\make@bigtimes{#1}}}}
\newcommand{\make@bigtimes}[1]{%
  \sbox\z@{$\m@th#1\sum$}%
  \setlength{\unitlength}{\wd\z@}%
  \begin{picture}(1,1)
    \roundcap
    \linethickness{.17ex}
    \Line(0.01,0)(.85,.99)
    \Line(0.01,.99)(.85,0.0)
  \end{picture}%
}
\makeatother

\newcommand{\Cpp}{C{\tt ++}}
\newcommand{\Alpaqa}{\textsc{Alpaqa}}
\newcommand{\alpaqa}{\textsc{alpaqa}}
\newcommand{\panoc}{\textsc{panoc}}
\newcommand{\Panoc}{\textsc{Panoc}}
\newcommand{\cutest}{\textsc{cute}st}
\newcommand{\casadi}{Cas\textsc{AD}i}
\newcommand{\ipopt}{\textsc{Ipopt}}
\newcommand{\qpalm}{\textsc{qpalm}}
\newcommand{\Qpalm}{\textsc{Qpalm}}

\author{\authorblockN{Pieter Pas, \and Mathijs Schuurmans, \and Panagiotis Patrinos}
    \thanks{The authors are with the Department
    of Electrical Engineering (\textsc{esat-stadius}), KU Leuven, 
    Kasteelpark Arenberg 10, 3001 Leuven, Belgium.
    Email: \texttt{\{pieter.pas, mathijs.schuurmans, panos.patrinos\}@esat.kuleuven.be}}%
    \thanks{This work was supported by: FWO projects: No. G0A0920N, No. G086318N; No. G086518N; 
    Fonds de la Recherche Scientifique -- FNRS, the Fonds Wetenschappelijk Onderzoek--Vlaanderen under EOS Project No. G0F6718N (SeLMA);
    Research Council KU Leuven C1 project No. C14/18/068; 
    European Union's Horizon 2020 research and innovation programme under the Marie Skłodowska-Curie grant agreement No. 953348;
    and Ford--KU Leuven Research Alliance project No. KUL0075.}
}%
\title{\Alpaqa{}: A matrix-free solver for nonlinear MPC and large-scale nonconvex optimization}

\begin{document}
\maketitle
\begin{abstract}
This paper presents \alpaqa{}, an open-source
\Cpp{} implementation of an augmented Lagrangian method for
nonconvex constrained numerical optimization, using the first-
order \panoc{} algorithm as inner solver. The implementation
is packaged as an easy-to-use library that can be used in \Cpp{}
and Python. Furthermore, two improvements to the \panoc{}
algorithm are proposed and their effectiveness is demonstrated
in NMPC applications and on the \cutest{} benchmarks for
numerical optimization.
The source code of the \alpaqa{} library is available at 
\href{https://github.com/kul-optec/alpaqa}{https://github.com/kul-optec/alpaqa}
and binary packages can be installed from 
\href{https://pypi.org/project/alpaqa}{https://pypi.org/project/alpaqa}.
\end{abstract}
\vspace{-1em}
\section{Introduction}

\subsection{Background and motivation}
Large-scale optimization problems arise naturally in many areas of science and engineering. As the models used 
in these different application domains are becoming increasingly complex, 
the need arises for solvers that can handle large, nonconvex optimization problems.
For instance, in the field of systems and control, the demand for improved autonomy in dynamical, complex 
and uncertain environments have boosted the importance of optimal control under nonlinear dynamics and long 
prediction horizons. This in turn yields nonconvex and large-scale optimization problems.

One of the core motivations for the present work is the requirement for real-time solution of nonconvex 
optimization problems in the context of nonlinear model predictive control (NMPC).
The real-time requirement, and often limited hardware capabilities available on embedded devices, result in 
a need for efficient solvers, both in terms of computation and memory.

Classical techniques for solving NMPC problems include sequential quadratic 
programming (SQP) and interior point (IP) methods. An advantage of SQP is that 
the solver can be warm-started using (a shifted version of) the solution of the 
previous time step, and techniques like real-time iteration (RTI)
require the solution of just a single QP at each time step. Interior point 
methods are popular general-purpose solvers, and state-of-the-art implementations such as 
\ipopt{} \cite{ipopt} are fast and easy to use. A disadvantage of interior point methods for model predictive control is that they cannot easily be warm-started, 
although some techniques do exist \cite{BensonHandeY2008Imfn}.

Another family of algorithms for constrained optimization that saw increased interest in the past decade are
Augmented Lagrangian methods (ALM) \cite{practicalaugmentedlagrangian}.
These methods involve the successive minimization of the augmented Lagrangian function, an 
exact penalty function for the original constrained optimization problem. The minimization problems in the inner loop of the ALM algorithm 
have simpler constraints than the 
original problem (often box constraints), and can be solved using a wide range of existing solvers. 
Thanks to the use of a penalty function, ALM can cope with large numbers of constraints. Moreover, by repeatedly solving similar problems, ALM naturally takes advantage of warm-starting the inner solver.
These properties render it particularly suitable for MPC applications, where a good initial guess is typically available.

In this work, ALM is used to tackle general nonconvex constraints, 
and the recently developed method called \panoc{} \textit{(Proximal Averaged Newton-type method for Optimal Control)} \cite{panoc}
is used as an inner solver. 
\Panoc{} is a quasi-Newton-accelerated first-order method, 
which -- despite its simple, inexpensive iterations -- 
enjoys locally superlinear convergence, without requiring the manipulation of large
Jacobian or Hessian matrices. 
Applying ALM and penalty-based methods combined with \panoc{} has recently proven successful within several applications in NMPC \cite{8550253,small_AerialNavigationObstructed_2019} and the Optimization Engine (\textsc{OpEn}) software \cite{open}, which shares some of its algorithmic foundations with the present work, although it is tailored more heavily towards embedded systems applications.

\subsection{Contributions}
The main contribution of this paper is \alpaqa{}
(\textit{Augmented Lagrangian Proximal Averaged Quasi-Newton Algorithms}), an 
open-source software package that provides efficient \Cpp{} implementations 
of an augmented Lagrangian method and different variants of \panoc{}. 
The solvers can be used both directly from \Cpp{} and from Python through a user-friendly interface based on \casadi{} \cite{Andersson2019}.

Furthermore, we propose two modifications to 
\panoc{} that significantly improve the practical performance and 
robustness of the algorithm: \Cref{sec:strucpanoc} exploits the structure of the 
derivatives of the objective function when applied to box-constrained problems, 
and \Cref{sec:impr-ls} introduces an alternative line search condition to help 
reject poor quasi-Newton steps that harm the practical convergence. 
In \Cref{sec:results}, the effectiveness of these modifications is demonstrated by 
applying the \alpaqa{} implementation to an NMPC problem and 
a subset of the \cutest{} benchmark collection.
\section*{Notation}
Let \([a,b]\) denote the closed interval from \(a\) to \(b\). \(\N_{[i,j]} \triangleq [i,j] \cap \N\) is the inclusive range of natural numbers from \(i\) to \(j\).
The set of extended real values is denoted as \(\overline \Re \triangleq \Re \cup \{+\infty\}\).
Uppercase calligraphic letters are used for index sets.
$x_i$ refers to the $i$'th component of $x \in \Re^n$. 
Given an index set $\mathcal{I} = \{i_1,\, \dots,\, i_m\} \subseteq \N_{[1,n]}$, we use the shorthand $x_{\mathcal{I}} = (x_{i_1},\, \dots,\, x_{i_m})$. For \(u,v \in \Re^n\), let \(u \le v\) denote the component-wise comparison.
Given a positive definite matrix \(\Sigma\), define the \(\Sigma\)-norm as
\(\|x\|_\Sigma \triangleq \sqrt{x^{\!\top} \Sigma\, x}\), 
and the distance in \(\Sigma\)-norm of a point \(x\) to a set \(C\) as 
\(\dist_\Sigma(x, C) \triangleq \min \left\{ \|x - v\|_\Sigma \;\mid\; v \in C\right\}\).
The proximal operator of function \(h : \Re^n \rightarrow \overline \Re\) is defined as \(\prox_{\gamma h}(x) \triangleq 
\min_w \big\{ h(w) + \tfrac{1}{2\gamma} \left\| w - x \right\|^2 \big\}\),
and the Moreau envelope
of \(h\) is defined as the value function of \(\prox_{\gamma h}\),
\(
    h^\gamma(x) \triangleq \inf_w \big\{ h(w) + \tfrac{1}{2\gamma} \left\| w - x \right\|^2 \big\}
\)
\cite[Def. 1.22]{RockafellarVariationalAnalysis}.

\section{Problem statement and preliminaries}
We propose a numerical
optimization solver for general nonconvex nonlinear programs of the form
\begin{equation}
    \begin{aligned}
        & \underset{x}{\text{minimize}}
        & & f(x) \\
        & \text{subject to}
        & & \underline{x} \le x \le \overline{x} \\
        &&& \underline{z} \le g(x) \le \overline{z},
    \end{aligned}
    \label{eq:problem-orig}
    \tag{P}
\end{equation}
with the objective function \(f : \Re^n \rightarrow \Re\)
and the constraints function \(g : \Re^n \rightarrow \Re^m\) 
possibly nonconvex.

Define the rectangular sets \(C\) and \(D\) as the cartesian products of 
one-dimensional closed intervals,
\(C_i \triangleq [\underline x_i, \overline x_i],\; C \triangleq 
C_1 \times \dots \times C_n\) and \(D_i \triangleq [\underline z_i, \overline z_i],\; D \triangleq 
D_1 \times \dots \times D_m\), 
such that the constraints of Problem \eqref{eq:problem-orig} can be written as
\(x \in C\) and \(g(x) \in D\).


For completeness, we briefly review the augmented Lagrangian method and the \panoc{} algorithm \cite{practicalaugmentedlagrangian,panoc}, including a minor modification to the ALM penalty factor, following \cite{qpalm}.
The methods described in this section serve as the algorithmic foundation for \alpaqa. In \Cref{sec:strucpanoc,sec:impr-ls}, we subsequently 
describe included improvements to these basic methods. 

\subsection{Augmented Lagrangian method}
The general constraints \(g(x) \in D\) in Problem \eqref{eq:problem-orig} 
are relaxed by introducing a slack variable \(z\) and by applying an augmented 
Lagrangian method to the following problem:
\begin{equation}
    \begin{aligned}
        & \underset{x \in C,\; z \in D}{\text{minimize}}
        & & f(x) \\
        & 
        \text{subject to}
        & & z = g(x).
    \end{aligned}
    \label{eq:problem-alm}
    \tag{P-ALM}
\end{equation}

Given a diagonal, positive definite \(m\)-by-\(m\) matrix \(\Sigma\),
define the augmented Lagrangian function
\(\Lagr_\Sigma : \reals^n \times \reals^m \times \reals^m \rightarrow \reals\)
with penalty factor \(\Sigma\) as
\begin{equation}
    \begin{aligned}
        \Lagr_\Sigma(x, z, y) \,\triangleq\, \Lagr(x, z, y) + \tfrac{1}{2} \left\|g(x) - z\right\|^2_\Sigma,
    \end{aligned}
\end{equation}
where $\Lagr(x,z,y) = f(x) +  \langle y, g(x) - z\rangle$ denotes the standard Lagrangian. 

The augmented Lagrangian method for solving Problem \eqref{eq:problem-alm}
then consists of repeatedly
\begin{inlinelist*}
\item \label{alm:step1} minimizing \(\Lagr_\Sigma\) with respect to the decision variables \(x\) and the slack variables \(z\)
\item \label{alm:step2} updating the Lagrange multipliers \(y\)
\item \label{alm:step3} increasing the penalty factors \(\Sigma_{ii}\) corresponding to constraints with high violation.
This procedure is shown in \Cref{alg:alm}. 
The minimization problem in 
step \ref{alm:step1} is solved up to a tolerance \(\varepsilon^\nu\) with \(\lim_{\nu\rightarrow\infty} \varepsilon^\nu = 0\) 
such that
\(\left\|x^\nu - \Pi_C\big(x^\nu - \nabla_{\!x}\! \Lagr_\Sigma(x^\nu, z^\nu, y^\nu)\big)\right\|_\infty \le \varepsilon^\nu\), with \(x^\nu \in C,\,z^\nu \in D\).
The solution \((x^\nu,\,y^\nu)\)
is returned when the termination criteria \(\left\|g(x^\nu) - z^\nu\right\|_\infty \le \delta\) and \(\varepsilon^\nu \le \varepsilon\) are satisfied, for given tolerances \(\varepsilon\) and \(\delta\).
\end{inlinelist*}

\begin{algorithm}
    \caption{Augmented Lagrangian Method}
    \label{alg:alm}
    \algsize
    \begin{algorithmic}[0]\vspace{0.66em}
        \Procedure {ALM}{$x^0$, $y^0$, $\Sigma^0$}\vspace{0.66em}
        \For{$\nu = 1, 2, \dots$}\vspace{0.66em}
        \State $(x^\nu, z^\nu) \leftarrow \underset{x\in C,\; z\in D}{\text{\textbf{argmin}}}\quad \Lagr_{\Sigma^{\nu-1}}(x,\, z;\, y^{\nu-1})$ \Comment{\ref{alm:step1}}
        \State $y^\nu \leftarrow \Pi_Y\left( y^{\nu-1} + \Sigma^{\nu-1} \left(g(x^\nu) - z^\nu\right)\right)$ \Comment{\ref{alm:step2}}\vspace{0.66em}
        \State $\Sigma^\nu \leftarrow$ \textsc{\footnotesize Update\_\(\Sigma\)}$(x^\nu,\, z^\nu,\, x^{\nu-1},\, z^{\nu-1},\, \Sigma^{\nu-1}$)\Comment{\ref{alm:step3}}
        \EndFor
        \EndProcedure
        \Statex
        \Function {Update\_\(\Sigma\)}{$x^\nu$, $z^\nu$, $x^{\nu-1}$, $z^{\nu-1}$, $\Sigma^{\nu-1}$}
        \State $e^\nu \leftarrow g(x^\nu) - z^\nu$
        \State $e^{\nu-1} \leftarrow g(x^{\nu-1}) - z^{\nu-1}$
        \State $\Sigma^\nu \leftarrow \Sigma^{\nu-1}$
        \For{$i = 1, 2, \dots, m$}
        \If{$|e^\nu_i| > \theta |e^{\nu-1}_i|$}
        \State $\Sigma^\nu_{ii} \leftarrow \Sigma^{\nu-1}_{ii}\,\max \left\{1,\; \Delta \dfrac{|e^\nu_i|}{\|e^\nu\|_\infty}\right\}$
        \EndIf
        \EndFor
        \State \textbf{return} $\Sigma^\nu$
        \EndFunction
    \end{algorithmic}
\end{algorithm}

The use of a diagonal matrix \(\Sigma\) as penalty factor --
with a separate penalty factor for each constraint, rather 
than a single scalar for all constraints -- was inspired by the \qpalm{} solver 
\cite{qpalm}. \Qpalm{} uses the same heuristic to update the penalty factors, 
where the increase in the penalty for a given constraint depends on the 
violation of that constraint relative to the total constraint violation.

\subsection{\Panoc}

The inner minimization problem in step \ref{alm:step1} of \Cref{alg:alm} will be solved
by the \panoc{} algorithm \cite{panoc}. \Panoc{} solves problems of the form:
\begin{equation}
    \begin{aligned}
        & \underset{x}{\text{minimize}}
        & & \psi(x) + h(x),
    \end{aligned}
    \label{eq:problem-panoc}
    \tag{P-PANOC}
\end{equation}
where \(\psi : \Re^n \rightarrow \Re\) is real-valued and 
\(h : \Re^n \rightarrow \overline \Re\) allows efficient computation of the 
proximal operator.
Recall that in \ref{alm:step1}, our aim is to solve
\begin{equation}
    \begin{aligned}
        \min_{x\in C,\, z\in D} \Lagr_\Sigma(x, z, y)
        &= \min_{x\in C}\big\{  \psi_{\Sigma}(x;y)  \big\}
        -\tfrac{1}{2} \lVert y \rVert_{\Sigma^{-1}}^2,\hspace{-0.1em}
    \end{aligned}
\end{equation}
\vspace{-0.5em}
where
\vspace{0.2em}
\begin{equation}
    \begin{aligned}
        \psi_{\Sigma}(x;y) &\triangleq f(x)
        + \min_{z\in D} \left\{
        \tfrac{1}{2} \left\Vert z - \left(g(x) + \Sigma^{-1} y\right)\right\Vert_\Sigma^2
        \right\} \\
        &= f(x)
        + \tfrac{1}{2} \text{dist}_\Sigma^2 \left(
        g(x) + \Sigma^{-1}y, \ D
        \right).
    \end{aligned}
\end{equation}
The minimizer \(z^\nu = \Pi_D\left(g(x^\nu) + \Sigma^{-1} y\right)\) can 
be computed directly.
Setting \(\psi(x) = \psi_\Sigma(x; y)\) and \(h = \delta_C\), 
the indicator function of $C$, it is clear that \panoc{} is indeed applicable to
the problem at hand.

The main idea behind \panoc{} is to select the next iterate as a convex combination of a proximal gradient 
step and a quasi-Newton step: Iteration using proximal gradient steps is slow, but convergence to a stationary point
can be guaranteed, whereas the quasi-Newton steps are faster, but only convergent close to a local minimum.
A line search is used as a globalization strategy, preserving the global convergence of the proximal gradient 
algorithm, while preferring quasi-Newton steps whenever possible to speed up the algorithm. A popular quasi-Newton method is limited-memory BFGS (L--BFGS)
\cite[Algorithm 2]{nocedal-lbfgs}.

To set the stage for the proposed improvements in the following 
sections, we now introduce the necessary details regarding \panoc{}.
Let the forward-backward operator \(T_\gamma : \reals^n \rightrightarrows \reals^n\) denote the set of
proximal gradient steps in \(x\), i.e.,
\begin{equation}
    \begin{aligned}
        T_\gamma(x) &\triangleq \prox_{\gamma h}\big(x - \gamma\nabla\psi(x)\big) \\
        &= \argmin_w \left\{ h(w) + \tfrac{1}{2 \gamma}\left\|x - \gamma\nabla\psi(x) - w\right\|^2\right\}.
    \end{aligned}
\end{equation}
Local solutions of \eqref{eq:problem-panoc} are fixed points of \(T_\gamma\), 
and roots of the fixed-point residual \(R_\gamma(x) \triangleq \frac{1}{\gamma} \left(x - T_\gamma(x)\right)\).
Let \(H_k\) be an operator that is in some sense close to the inverse Jacobian
of the fixed-point residual \(R_\gamma(x)\)\footnote{Specifically, if \(H_k\) satisfies the Dennis-Mor\'e condition \cite{dennismore}, local superlinear convergence can be achieved. \cite[Thm.~III.5]{panoc}}.
Then for a given step size \(\gamma\), the next \panoc{} iterate is computed as
\begin{subequations}
    \begin{align}
        \hat x^k &\in T_\gamma(x^{k}) \\ 
        p^k &= \hat x^k - x^k \label{eq:def-step-p} \\
        q^k &= -H_k R_\gamma(x^k) = -\gamma^{-1} H_k p^k \\
        x^{k+1} &= x^k + (1-\tau)\, p^k + \tau\, q^k.
    \end{align}
\end{subequations}
The parameter \(\tau \in [0, 1]\) is selected using a backtracking line search over the forward-backward envelope (FBE) \cite{PatrinosPanagiotis2013PNmf}, an exact, continuous, 
real-valued penalty function for the inner problem \eqref{eq:problem-panoc}.
More specifically,
let the parameter
\(
\sigma \in \left( 0,\; \tfrac{1-\gamma L}{2\gamma} \right),
\) where the step size \(\gamma\) and the parameter \(L \in \left(0, \gamma^{-1}\right)\) are chosen such that the
following quadratic upper bound is satisfied:
\begin{equation} \label{eq:qub}
    \psi(\hat{x}^k) \le \psi(x^k) + \nabla\psi(x^k)^\top (\hat{x}^k - x^k) + \tfrac{L}{2} \lVert \hat{x}^k - x^k \rVert^2.
\end{equation}
Then, \panoc{} uses the sufficient decrease condition
\begin{equation} \label{eq:panoc-ls}
    \varphi_{\gamma}(x^{k+1}) \le \varphi_{\gamma}(x^k) - \sigma \lVert p^k \rVert^2, 
\end{equation}
where \(\varphi_\gamma : \reals^n \rightarrow \reals\) is the FBE 
defined as 
\begin{equation}\label{eq:fbe}
    \begin{aligned}
        \varphi_\gamma(x) &\triangleq \psi(x) - \tfrac{\gamma}{2} \| \nabla\psi(x) \|^2 + h^\gamma(x - \gamma \nabla\psi(x)),
    \end{aligned}
\end{equation}
and \(p^k\) is the proximal gradient step as defined in \eqref{eq:def-step-p}.

When these conditions hold, it can be shown that when \(\tau = 0\) or 
\(x^{k+1} = \hat x^k\),
the line search condition \eqref{eq:panoc-ls} is satisfied, and  
global convergence of the algorithm can be demonstrated by telescoping 
equation \eqref{eq:panoc-ls} \cite[Thm. III.3]{panoc}.

At each iterate, condition \eqref{eq:qub} is satisfied using the procedure 
outlined in \Cref{alg:step-size}. If the gradient of \(\psi\) is 
\(L_{\psi}\)-Lipschitz continuous, then the \texttt{while} loop will terminate, because 
condition \eqref{eq:qub} holds whenever \(L \ge L_{\psi}\). If \(\nabla \psi\)
is only locally Lipschitz and if the domain is bounded, then the maximum value of 
\(L\) can be bounded as well, again implying the termination of \Cref{alg:step-size}. \cite[Remark III.4]{panoc}
The original \panoc{} algorithm is listed in \Cref{alg:old-panoc}.

\begin{algorithm}
    \caption{Step size selection}
    \label{alg:step-size}
    \algsize
    \begin{algorithmic}[0]\vspace{0.66em}
        \Function {update\_step\_size}{$x$, $\gamma$, $L$}
        \State $\hat x \leftarrow T_\gamma(x)$, \hspace{0.5em} $p \leftarrow \hat x - x$
        \While{
            $\psi(\hat{x}) > \psi(x) + \nabla\psi(x)^\top p + \tfrac{L}{2} \lVert p \rVert^2$
        }
            \State $\gamma \leftarrow \gamma / 2$, \hspace{0.5em} $L \leftarrow 2L$
            \State $\hat x \leftarrow T_\gamma(x)$, \hspace{0.5em} $p \leftarrow \hat x - x$
            \EndWhile
        \State \textbf{return} $\gamma$, $L$
        \EndFunction
    \end{algorithmic}
\end{algorithm}

\algdef{SE}[DOWHILE]{Do}{doWhile}{\algorithmicdo}[1]{\algorithmicwhile\ #1}%

\begin{algorithm}
    \caption{\panoc}
    \label{alg:old-panoc}
    \algsize
    \begin{algorithmic}[1]\vspace{0.66em}
    \Procedure{\panoc{}}{$x^0$, $\gamma_0$, $L_0$}
    \vspace{0.66em}
    \For{\(k \leftarrow 0, 1, 2, \dots\)}\vspace{0.66em}
        \State \(p^k \leftarrow T_{\gamma_k}(x^k) - x^k\) \Comment{Proximal gradient step} 
        \State \label{step:quasi-newt}\(q^k \leftarrow -\gamma^{-1}H_k p^k\) \Comment{Quasi-Newton step}\vspace{0.66em}
        \State \(\tau \leftarrow 1\)
        \Do \Comment{Line search}
            \State \(x^{k+1} \leftarrow x^k + (1-\tau)\, p^k + \tau\, q^k\)
            \State \(\tau \leftarrow \tau / 2\)
        \doWhile{\(\varphi_{\boldsymbol{\gamma_{k}}}(x^{k+1}) > \varphi_{\gamma_k}(x^k) - \sigma_k \lVert p^k
            \rVert^2\)}
        \vspace{0.66em}
        \State $\gamma_{k+1}$, $L_{k+1}$ $\leftarrow$ \textsc{{\footnotesize update\_step\_size}($x^{k+1}$, $\gamma_{k}$, $L_{k}$)}
    \EndFor
    \EndProcedure
    \end{algorithmic}
\end{algorithm}
\section{Structured \panoc}\label{sec:strucpanoc}

In the original \panoc{} algorithm, $h$ was essentially allowed to be any proximable mapping. Since for our purposes, $h$ in \eqref{eq:problem-panoc} is taken to be the indicator of $C$, 
its proximal operator reduces to a projection onto the box, i.e.,  
\(\prox_{\gamma \delta_C}(x)
 = \argmin_{w\in C} \big\{ \tfrac{1}{2 \gamma}\left\|w - x\right\|^2\big\}
= \Pi_C(x)\). This results in an inexpensive component-wise operation and furthermore induces some additional structure, which we may exploit to
replace line~\ref{step:quasi-newt} in \Cref{alg:old-panoc}, which was originally implemented using a general L--BFGS strategy. 
This results in a more efficient variant of \panoc{}, which is better tailored to our particular purpose.

To this end let us first define the set of indices corresponding to active box constraints on $x$ after a gradient step, 
\begin{equation} 
    \begin{aligned}
        \sK(x) &\triangleq \left\{\; i \in \N_{[1,n]} \;\left|\; \begin{aligned} &x_i - \gamma \nabla_{x_i} \psi(x) \le \underline x_i \\
        \;\;\vee\;\; &\overline x_i \le x_i - \gamma\nabla_{x_i}\psi(x) \end{aligned} \right. \right\}\hspace{-0.5em}
    \end{aligned}
\end{equation}
and its complement \(\sJ(x) \triangleq \N_{[1,n]} \setminus \sK(x)\) (inactive constraint indices). For ease of notation, let \(P_{\sK\sJ}\) denote a
row permutation matrix that reorders all rows with indices \(k\in \sK(x)\) before the rows with indices \(j \in \sJ(x)\).

A straightforward calculation reveals that 
if \(\nabla \psi\) and the fixed-point residual \(R_\gamma = \gamma^{-1} (I - T_{\gamma})\) are differentiable at \(x\), the Jacobian \(\jac R_{\gamma}(x)\) is given by
\begin{equation}\label{eq:lem-jac-fpr}
    \jac R_\gamma(x) = P_{\sK \sJ}^\top
    \begin{pmatrix}
        \gamma^{-1} I & 0 \\
        \nabla^2_{\!x_\sJ x_\sK}\! \psi(x) & \nabla^2_{\!x_\sJ x_\sJ}\! \psi(x)
    \end{pmatrix}
    P_{\sK \sJ},
\end{equation}
using which we may apply a semismooth Newton method to the problem of finding roots of \(R_\gamma\),
replacing the quasi-Newton step \(q^k = -H_k R_\gamma(x^k)\).
The Newton step \(q\) is computed as the solution to the system of equations 
\begin{equation}
    \begin{pmatrix}
        \gamma^{-1} I & 0 \\
        \nabla^2_{\!x_\sJ x_\sK}\! \psi(x) & \nabla^2_{\!x_\sJ x_\sJ}\! \psi(x)
    \end{pmatrix}
    \;
    P_{\sK\sJ}\, q
    = - P_{\sK\sJ}\, R_\gamma(x).
\end{equation}
Defining \(\begin{pmatrix}  q_\sK \\ q_\sJ \end{pmatrix} \triangleq P_{\sK\sJ}\, q\), this can be written as 
\begin{subequations}
    \begin{align}
        \label{eq:61ab}
        q_\sK &= x_\sK - T_{\gamma}(x)_\sK \\
        \label{eq:61c}
        \nabla^2_{\!x_\sJ x_\sJ}\! \psi(x)\, q_\sJ &= -\nabla_{\!x_\sJ}\psi(x) - \nabla^2_{\!x_\sJ x_\sK}\! \psi(x)\, q_\sK.
    \end{align}
\end{subequations}
\noindent%
Here, the Hessian of \(\psi\) is given by
\begin{equation}
    \begin{aligned}\label{eq:hes-obj-1}
        \nabla^2 \psi(x) =\;&
        \nabla^2 f(x) + \sum_{i\in\mathcal{A}(x)} \nabla^2 g_i(x)\, \hat y_i(x) \\
        &\phantom{\nabla^2 f(x)} + \sum_{i\in\mathcal{A}(x)} \Sigma_{ii}\,\nabla g_i(x) \nabla g_i(x)^\top
        \end{aligned}
        \vspace{-1em}
\end{equation}
where 
\begin{equation} \label{eq:hess-active-set-I}
    \begin{aligned}
    \mathcal{A}(x) = \left\{\; i \in \N_{[1,m]} \;\left|\; \begin{aligned} &g_i(x) + \Sigma_{ii}^{-1} y_i < \underline z_i \\
    \;\;\vee\;\; &\overline z_i < g_i(x) + \Sigma_{ii}^{-1} y_i \end{aligned} \right.\right\}
    \end{aligned}
\end{equation}
is the set of indices of the active (nonlinear) constraints, and \\
\(
    \hat y(x) \triangleq \Sigma\, \left(g(x) + \Sigma^{-1}y - \Pi_D\left(g(x) + \Sigma^{-1}y\right)\right).
\)

Equation \eqref{eq:61ab} can be computed without any additional work, because these components 
of the Newton step are equal to the projected gradient step,
\(q_\sK = p_\sK\). \\
Equation \eqref{eq:61c} on the other hand requires the solution of a (smaller) linear 
system using part of the Hessian of the augmented Lagrangian function, \(\nabla^2_{\!x_\sJ x_\sJ}\!\psi\).

\Cref{subsec:aughess-newton,subsec:lagrhess-newton} present different strategies for solving the system \eqref{eq:61c}.
\subsection{Direct solution using the Hessian of the augmented Lagrangian}
\label{subsec:aughess-newton}
A first approach to compute the step \(q_\sJ\) in equation~\eqref{eq:61c} is to
solve it directly using standard techniques such as \(LDL^\top\) or Cholesky 
factorization. If the block
\(\nabla^2_{\!x_\sJ x_\sJ}\!\psi(x)\) of the Hessian
is not positive definite, a modified Cholesky factorization 
\cite{FangHaw-ren2008MCaa}, or a conjugate gradients method \cite[§7]{nocedal_numopt} can be used.

A disadvantage of this approach is that, in general, the Hessian 
of the augmented Lagrangian \(\nabla^2 \psi\)
could be dense, even if the Hessian of the Lagrangian \(\hessxx \Lagr\) is sparse,
because of the rightmost term in equation \eqref{eq:hes-obj-1}.
For large-scale optimization problems, computing, storing and factorizing this dense Hessian is 
often prohibitively expensive.

\subsection{Solution using the Hessian of the Lagrangian and the Jacobian of the constraints}
\label{subsec:lagrhess-newton}
Using the equation~\eqref{eq:hes-obj-1}, the Hessian of the augmented 
Lagrangian can be written in terms of the Hessian of the Lagrangian and the 
Jacobian of the constraints, with the goal of exploiting the sparsity of these matrices.

For ease of notation, define the matrices
\newcommand{\Hjj}{\ensuremath{H_{\!\sJ}}}
\newcommand{\Gj}{\ensuremath{G_{\!\sJ}}}
\newcommand{\xj}{\ensuremath{{x_{\!\sJ}}}}
\newcommand{\xk}{\ensuremath{{x_{\!\sK}}}}
\ifArxiv
\begin{equation*}
\else
\\ \(
\fi
    \Hjj \triangleq \Big. \nabla_{\!x_\sJ x_\sJ}^2\! \Lagr(x, y) \Big|_{\big(x,\, \hat y(x)\big)} \;\text{and}\;\;
    \Gj \triangleq \nabla_{\!\xj} g_\mathcal{A}(x).
\ifArxiv
\end{equation*}
\else
\) \\
\fi
\ifArxiv
\else
Then \eqref{eq:61c} is equivalent to the following linear system
\fi
\begin{equation*}
    \begin{pmatrix}
        \Hjj & \Gj \\ \Gj^\top & -\Sigma_{\!\mathcal{A}}^{-1}
    \end{pmatrix} 
    \begin{pmatrix}
        q_\sJ \\ v
    \end{pmatrix}
    = -
    \begin{pmatrix}
        \nabla_\xj\! f(x) + \nabla^2_{\!\xj \xk}\!\psi(x)\, q_\sK \\
        \Sigma_{\!\mathcal{A}}^{-1}\, \hat y_\mathcal{A}(x)
    \end{pmatrix}
    \ifArxiv\else , \fi
\end{equation*}
\ifArxiv
\else
where \(v = \Sigma_{\!\mathcal{A}} \, \Gj^\top\, q_\sJ + \hat y_\mathcal{A}(x)\).
\fi
This system is larger than \eqref{eq:61c}, but its sparsity pattern consists of the sparsity of the Hessian of the 
Lagrangian and the Jacobian of the constraints. It can be solved using symmetric indefinite (sparse) linear
solvers, such as \texttt{MA57} \cite{DuffIain2004Mcft}.

\subsection{Approximate solution using quasi-Newton methods}\label{subsec:struclbfgs}

In the spirit of the original \panoc{} algorithm, one could use a quasi-Newton method such as L--BFGS to approximate the inverse of
\(\nabla^2_{\!x_\sJ x_\sJ}\!\psi(x)\) instead of explicitly solving a linear system.
Since the set of active constraints \(\sK(x)\) may change at any iteration, the L--BFGS algorithm \cite[Algorithm 7.4]{nocedal_numopt}
has to be updated to account for this:
The full vectors \(s^k \triangleq x^{k+1} - x^k \in \Re^n\) and \(y^k = \nabla\psi(x^{k+1}) - \nabla\psi(x^k) \in \Re^n\)
are stored,
and when the L--BFGS estimate is applied to a vector, only the components of \(s^k\) and 
\(y^k\) with indices in the index set \(\sJ(x^k)\) are used. However, this means that the curvature condition 
\(y^{k\top}_{\!\sJ}\! s^k_{\!\sJ} > 0\) that ensures positive definiteness of the L--BFGS Hessian estimate \cite[Eq. 6.7]{nocedal_numopt}
cannot be verified when the vectors \(s^k\) and \(y^k\) are added, the condition 
has to be checked later when applying the estimate.

The last term of equation~\eqref{eq:61c}, the Hessian-vector product \(\nabla^2_{\!x_\sJ x_\sK}\!\psi(x)\, q_\sK\) can 
be computed using algorithmic differentiation, by computing finite differences of \(\nabla\psi\), or left out entirely.
By leaving out this term, the method becomes equivalent to applying L--BFGS to the problem of minimizing \(\psi\) without box 
constraints, but only applied to the components of the decision variable for which the box constraints of the original problem are inactive.\\

The strategies from \Cref{subsec:aughess-newton,subsec:lagrhess-newton} give rise to 
second-order methods and require specialized linear solvers. For this reason, 
\alpaqa{} implements the quasi-Newton method from \Cref{subsec:struclbfgs}.
\section{Improved line search}\label{sec:impr-ls}

A practical issue that occurs with the original \panoc{} algorithm is that 
when L--BFGS produces a step \(q^k\) of low quality (i.e., a step that results 
in a significant increase of the objective function and/or constraint violation),
it is 
sometimes accepted by the line search procedure. The reason for this counterintuitive 
acceptance is that the quadratic upper bound of equation~\eqref{eq:qub} 
is not necessarily satisfied between at \(x^{k+1}\!,\) and as a result, 
the FBE 
\(\varphi_\gamma\) cannot be bounded from below \cite[Proposition 2.1]{ThemelisAndreas2018Feft}. This means that the FBE might 
decrease significantly by accepting the step \(q^k\!,\) even when the original objective
\(\psi + h\) increases. 
This phenomenon does not affect the global convergence of the algorithm. However, in the interest of practical performance,
it is beneficial to reject these low-quality L--BFGS steps, as two problems arise when such steps are accepted:
\begin{inlinelist*}
    \item it is likely that the next iterate
will be (much) farther away from the optimum than before
    \item the local curvature of \(\psi\) might be much larger at \(x^{k+1}\), demanding a 
significant reduction of the step size \(\gamma\) in the next iteration
in order to satisfy equation~\eqref{eq:qub}.
\end{inlinelist*}
The combined effect of 
ending up far away from the optimum \textit{and} having to advance with small steps is 
detrimental for performance.

\subsection{An improved, stricter line search condition}

The solution proposed here is to update the step size at the candidate iterate
\(x^{k+1} = x^k + (1-\tau)\, p^k + \tau\, q^k\)
using \Cref{alg:step-size}
\textit{before} evaluating the line search condition, and to use this new step 
size \(\gamma_{k+1}\) in the left-hand side of the line search condition
\eqref{eq:panoc-ls}, rather
than the old step size \(\gamma_k\).
As a consequence, candidate steps that would
result in a significant reduction of the step size are penalized by the new line 
search. 
If the candidate \(x^{k+1}\) does not require step size reduction,
\(\gamma_{k+1}\) is equal to \(\gamma_k\) and the modified algorithm is
equivalent to the original \panoc{} algorithm.

Specifically, the modified line search accepts the candidate iterate 
\(x^{k+1} = x^k + (1-\tau)\, p^k + \tau\, q^k\) if
\begin{equation} \label{eq:panoc-ls2}
    \varphi_{\gamma_{k+1}}(x^{k+1}) \le \varphi_{\gamma_{k}}(x^k) - \sigma \lVert p^k \rVert^2,
\end{equation}
where the step sizes \(\gamma\) and the parameters \(L \in \left(0, \gamma^{-1}\right)\) are chosen such that the
following quadratic upper bounds are satisfied in the current iterate \(x^k\)
and in the candidate iterate \(x^{k+1}\):
\begin{equation*} \label{eq:qub2}
    \begin{aligned}
        \psi(\hat{x}^k) &\le \psi(x^k) + \nabla\psi(x^k)^\top p^k + \tfrac{L_k}{2} \lVert p^k \rVert^2,\\
        \psi(\hat{x}^{k+1}) &\le \psi(x^{k+1}) + \nabla\psi(x^{k+1})^\top p^{k+1} + \tfrac{L_{k+1}}{2} \lVert p^{k+1} \rVert^2,
    \end{aligned}
\end{equation*}
with the usual definitions of
\(\hat{x}^k \triangleq T_{\gamma_k}(x^k)\), \(p^k \triangleq \hat{x}^k - x^k\) and 
\(\hat{x}^{k+1} \triangleq T_{\gamma_{k+1}}(x^{k+1})\), \(p^{k+1} \triangleq \hat{x}^{k+1} - x^{k+1}\).

The pseudocode for the \panoc{} algorithm with this modified line search is listed in \Cref{alg:new-panoc-detail};
compare the location of the call to \textsc{update\_step\_size} to the original version in \Cref{alg:old-panoc}.

\begin{algorithm}
    \caption{\panoc{} with improved line search condition}
    \label{alg:new-panoc-detail}
    \algsize
    \begin{algorithmic}[1]\vspace{0.66em}
    \Procedure{\panoc}{$x^0$, $\gamma_0$, $L_0$}
    \vspace{0.66em}
    \For{\(k \leftarrow 0, 1, 2, \dots\)}\vspace{0.66em}
        \State \(p^k \leftarrow T_{\gamma_k}(x^k) - x^k\) \Comment{Proximal gradient step} 
        \State \(q^k \leftarrow -\gamma^{-1}H_k p^k\) \Comment{Quasi-Newton step}\vspace{0.66em}
        \State \(\tau \leftarrow 1\)
        \Do \Comment{Line search}
            \If{\(\tau < \tau_\text{min}\)} 
                \State \(\tau \leftarrow 0\)
            \EndIf
            \State \(x^{k+1} \leftarrow x^k + (1-\tau)\, p^k + \tau\, q^k\)
            \State $\gamma_{k+1}$, $L_{k+1}$ $\leftarrow$ \textsc{{\footnotesize update\_step\_size}($x^{k+1}$, $\gamma_k$, $L_k$)}
            \State \(\tau \leftarrow \tau / 2\)
        \doWhile{\(\varphi_{\boldsymbol{\gamma_{k+1}}}(x^{k+1}) > \varphi_{\gamma_k}(x^k) - \sigma_k \lVert p^k
            \rVert^2\) and \(\tau \neq 0\)}
    \EndFor
    \EndProcedure
    \end{algorithmic}
\end{algorithm}

\subsection{Convergence of \panoc{} with improved line search}

The line search condition presented in the previous section is stricter than the original one because of the following trivial property 
of the FBE:
\begin{equation}\label{eq:ineq-gamma-fbe}
    \gamma_{k+1} \le \gamma_k \;\Rightarrow\; \varphi_{\gamma_{k+1}} \ge \varphi_{\gamma_k}.
\end{equation}
Since the new line search condition implies the original one, the global convergence results 
carry over to this improved algorithm.

One caveat is that proving termination of the line search loop is no longer 
possible, because it is not necessarily the case that 
\(\varphi_{\gamma_{k+1}}(\hat x^{k}) \le \varphi_{\gamma_{k}}(x^k) - \sigma \lVert p^k \rVert^2\).
Fortunately, this is not an issue in practice, since one could simply enforce termination 
after a finite number of line search iterations and accept the proximal gradient 
iterate \(\hat x^k\) regardless of whether it satisfies 
condition \eqref{eq:panoc-ls2}. Given that \(\hat x^k\) necessarily satisfies condition \eqref{eq:panoc-ls}, 
this does not affect the global convergence.
In fact, practical implementations of the original \panoc{} algorithm already limit the number
of line search iterations for 
performance reasons and to avoid numerical issues that might cause the line search not to 
terminate in practice.
For example, in \Cref{alg:new-panoc-detail}, 
this is achieved by setting \(\tau\) equal to zero when it 
reaches a certain threshold \(\tau_\text{min} \in (0, 1)\).

\section{Applications and experimental results}\label{sec:results}

In this section, the \alpaqa{} library is applied to different nonconvex
optimization problems, and a comparison is made of the performance and 
robustness of the proposed modifications.

\subsection{NMPC of a hanging chain}

\begin{figure}[ht!]
\centering 
\includegraphics[width=0.85\linewidth, clip,trim=0.7cm 0.6cm 0 0.2cm]{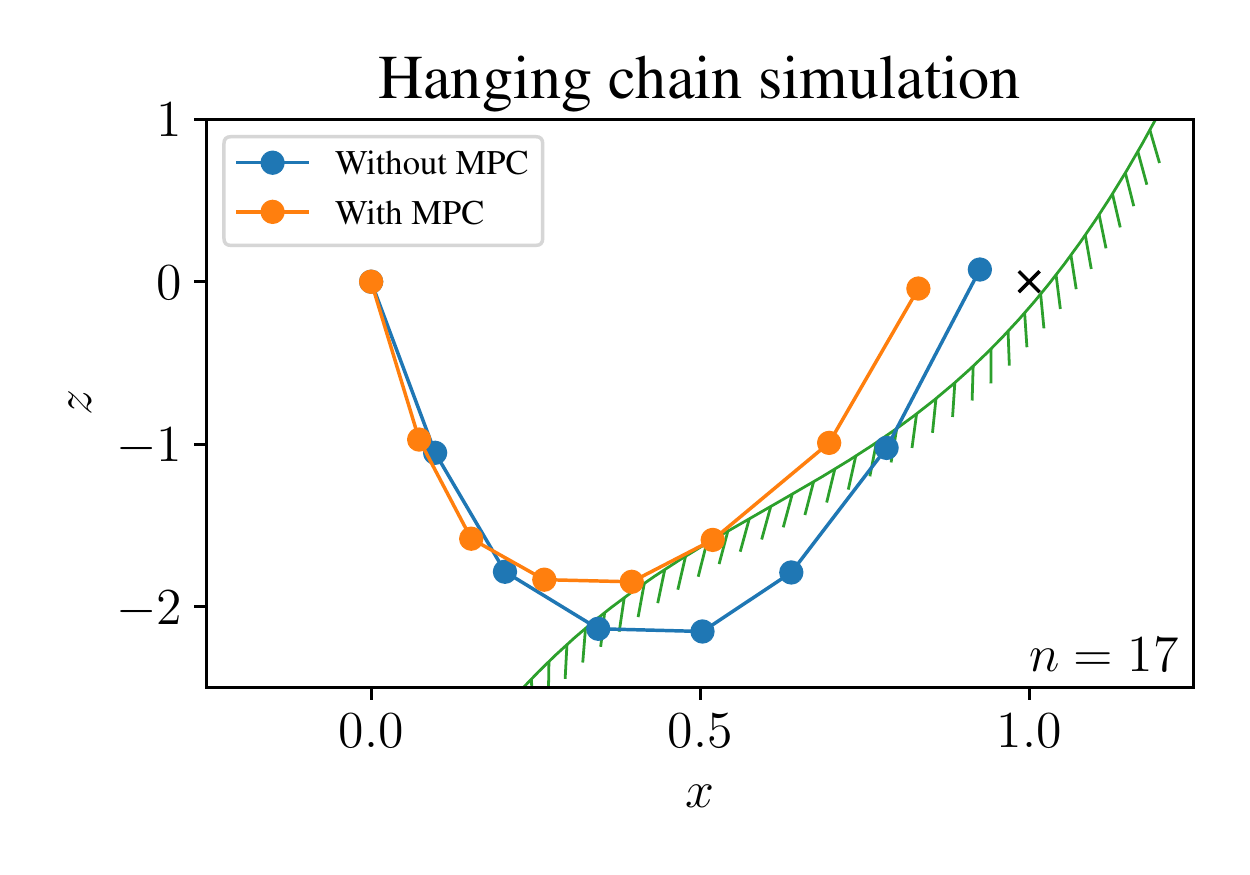}
\caption{Snapshot of the MPC simulation of the hanging chain at time step 17,
projected along the $y$-axis. The cubic state constraint is shown in green.}
\label{fig:chain}
\end{figure}

The hanging chain model from \cite{chainmodel} is used as a benchmark:
six balls are connected using seven springs. The first ball 
is connected to the origin using a spring, and the last ball is connected to 
an actuator using another spring. The origin remains fixed, the actuator can be 
moved through three-dimensional space with a maximum velocity of $1\,\mathrm{m}/\mathrm{s}$
in each of the three components. The model is discretized using a fourth-order 
Runge-Kutta method with a time step of $T_s = 0.05\,\mathrm{s}$. The springs have a 
spring constant $D = 1.6\,\mathrm{N/m}$, a rest length $L = 0.0055\,\mathrm{m}$
and the balls have a mass $m = 0.03\,\mathrm{kg}$. Additionally, 
the positions of the balls and the actuator are constrained by
\(x_z \ge c\,(x_x-a)^3 + d\,(x_x - a) + b\)
with \(a = 0.6, b = -1.4, c = 5, d = 2.2\).
After setting up the balls equidistantly between 0 and 1 on the $x$-axis, the
system is perturbed by an input of \(u = (-0.5,\,0.5,\,0.5)\) for three time 
steps.
An optimal control problem of horizon 40 is formulated, using 
\(L(x, u) = \alpha\,\|x_\text{actuator} - x_\text{end}\|_2^2 + 
\beta \sum_{i=1}^6 \| \dot x^{(i)} \|_2^2 + \gamma\,\|u\|_2^2\)
as the stage cost 
function,
with \(x_\text{actuator}\) the position vector of the actuator, \(x_\text{end}\)
the target position for the actuator, \(\dot x^{(i)}\) the velocity vector of 
the \(i\)-th ball, and \(u\) the velocity of the actuator, which is the 
input to the system. The goal of the MPC controller is to drive the system to 
the equilibrium state with the actuator at \(x_\text{end} = (1,\,0,\,0)\).
The specific weighting factors used are \(\alpha = 25\,\mathrm{s}^{-2}, \beta = 1\) and \(\gamma =
0.01\).
A single-shooting 
formulation is used, using ALM for the general state constraints and \panoc{} 
for the box constraints on the actuator velocity.
\cref{fig:chain} shows a schematic representation of the chain and the 
constraints.

\subsection{Improvements to \panoc{}}

A total of six different solver configurations are compared: three solvers 
-- the original \panoc{} solver as described in \cite{panoc}, 
structured \panoc{} using L--BFGS with inclusion of the Hessian-vector product
in equation \eqref{eq:61c}, and structured \panoc{} using L--BFGS without the 
rightmost term in \eqref{eq:61c} -- are all tested with the original line 
search condition and with the improved line search condition. 
The NMPC controller is simulated for 10 time steps, and the total 
number of function and gradient evaluations, the number of \panoc{} iterations, 
and the run time are reported.
A first experiment investigates the performance of the solvers when 
an initial guess of all zeros is used at each time step, whereas in a second 
experiment, the solver is warm-started using the solution and the Lagrange
multipliers from the previous time step, shifted by one time step. The results 
are shown in 
\ifMergefig
\cref{fig:cold-warm-start}.
\else
\cref{fig:cold-start,fig:warm-start} respectively.
\fi

\definecolor{pal1}{RGB}{16,40,108}
\definecolor{pal2}{RGB}{10,140,196}
\definecolor{pal3}{RGB}{0,121,44}
\definecolor{pal4}{RGB}{147,184,18}
\definecolor{pal5}{RGB}{226,6,25}
\definecolor{pal6}{RGB}{237,130,1}

\begin{figure}[bt!]
\vspace{0.6em}
\newlength{\myfigheight}
\setlength{\myfigheight}{0.45\linewidth}

\centering 
\pgfplotstableread[col sep=comma]{data/cold10-2/eval.csv}\table
    \begin{tikzpicture}
    \begin{axis}[
            SmallBarPlot,
            xticklabels from table={\table}{label},
            ylabel=Number of evaluations,
            xlabel={Function and gradient evaluations},
            title={\large Experiment I: Cold start},
            height=\myfigheight,
            legend style={at={(0.5,1.5), text=\scriptsize}, draw opacity=0.2, anchor=south},
            legend cell align=left,
        ]
        \addplot [fill=pal1]
            table [x expr=\coordindex, y=s0] {\table};
        \addlegendentry{\hspace{0.5em}Original \panoc{}}
        \addplot [fill=pal2]
            table [x expr=\coordindex, y=s1] {\table};
        \addlegendentry{\hspace{0.5em}Original \panoc{} with improved line search}
        \addplot [fill=pal3]
            table [x expr=\coordindex, y=s2] {\table};
        \addlegendentry{\hspace{0.5em}Structured \panoc{}}
        \addplot [fill=pal4]
            table [x expr=\coordindex, y=s3] {\table};
        \addlegendentry{\hspace{0.5em}Structured \panoc{} with improved line search}
        \addplot [fill=pal5]
            table [x expr=\coordindex, y=s4] {\table};
        \addlegendentry{\hspace{0.5em}Approx. structured \panoc{}}
        \addplot [fill=pal6]
            table [x expr=\coordindex, y=s5] {\table};
        \addlegendentry{\hspace{0.5em}Approx. structured \panoc{} with improved line search}
        \end{axis}
    \end{tikzpicture}
    \pgfplotstableread[col sep=comma]{data/cold10-2/inner_it.csv}\table
    \begin{tikzpicture}
        \begin{axis}[
            SmallBarPlot,
            xticklabels from table={\table}{label},
            ylabel=Number of iterations,
            xlabel=,
            width=0.45\linewidth,
            height=\myfigheight,
        ]
        \addplot [fill=pal1]
        table [x expr=\coordindex, y=s0] {\table};
        \addlegendentry{\hspace{0.5em}Original \panoc{}}
        \addplot [fill=pal2]
        table [x expr=\coordindex, y=s1] {\table};
        \addlegendentry{\hspace{0.5em}Original \panoc{} with improved line search}
        \addplot [fill=pal3]
        table [x expr=\coordindex, y=s2] {\table};
        \addlegendentry{\hspace{0.5em}Structured \panoc{}}
        \addplot [fill=pal4]
        table [x expr=\coordindex, y=s3] {\table};
        \addlegendentry{\hspace{0.5em}Structured \panoc{} with improved line search}
        \addplot [fill=pal5]
            table [x expr=\coordindex, y=s4] {\table};
        \addlegendentry{\hspace{0.5em}Approx. structured \panoc{}}
        \addplot [fill=pal6]
            table [x expr=\coordindex, y=s5] {\table};
        \addlegendentry{\hspace{0.5em}Approx. structured \panoc{} with improved line search}
        \legend{}
        \end{axis}
    \end{tikzpicture}
    \pgfplotstableread[col sep=comma]{data/cold10-2/runtime.csv}\table
    \begin{tikzpicture}
        \begin{axis}[
            SmallBarPlot,
            xticklabels from table={\table}{label},
            ylabel=Seconds,
            xlabel=,
            width=0.45\linewidth,
            height=\myfigheight,
        ]
        \addplot [fill=pal1]
            table [x expr=\coordindex, y=s0] {\table};
        \addlegendentry{\hspace{0.5em}Original \panoc}
        \addplot [fill=pal2]
            table [x expr=\coordindex, y=s1] {\table};
        \addlegendentry{\hspace{0.5em}Original \panoc{} with improved line search}
        \addplot [fill=pal3]
            table [x expr=\coordindex, y=s2] {\table};
        \addlegendentry{\hspace{0.5em}Structured \panoc}
        \addplot [fill=pal4]
            table [x expr=\coordindex, y=s3] {\table};
        \addlegendentry{\hspace{0.5em}Structured \panoc{} with improved line search}
        \addplot [fill=pal5]
            table [x expr=\coordindex, y=s4] {\table};
        \addlegendentry{\hspace{0.5em}Approx. structured \panoc}
        \addplot [fill=pal6]
            table [x expr=\coordindex, y=s5] {\table};
        \addlegendentry{\hspace{0.5em}Approx. structured \panoc{} with improved line search}
        \legend{}
        \end{axis}
    \end{tikzpicture}
    \ifMergefig
    \else
    \caption{Experimental results of the six solvers over 10 NMPC time steps 
    of the chain model, without warm-starting.}
    \label{fig:cold-start}
\end{figure}

\begin{figure}
    \fi
\pgfplotstableread[col sep=comma]{data/warm10-2/eval.csv}\table
    \begin{tikzpicture}
    \begin{axis}[
            SmallBarPlot,
            xticklabels from table={\table}{label},
            ylabel=Number of evaluations,
            xlabel={Function and gradient evaluations},
            title={\large Experiment II: Warm start},
            height=\myfigheight,
        ]
        \addplot [fill=pal1]
            table [x expr=\coordindex, y=s0] {\table};
        \addlegendentry{\hspace{0.5em}Original \panoc}
        \addplot [fill=pal2]
            table [x expr=\coordindex, y=s1] {\table};
        \addlegendentry{\hspace{0.5em}Original \panoc{} with improved line search}
        \addplot [fill=pal3]
            table [x expr=\coordindex, y=s2] {\table};
        \addlegendentry{\hspace{0.5em}Structured \panoc}
        \addplot [fill=pal4]
            table [x expr=\coordindex, y=s3] {\table};
        \addlegendentry{\hspace{0.5em}Structured \panoc{} with improved line search}
        \addplot [fill=pal5]
            table [x expr=\coordindex, y=s4] {\table};
        \addlegendentry{\hspace{0.5em}Approx. structured \panoc}
        \addplot [fill=pal6]
            table [x expr=\coordindex, y=s5] {\table};
        \addlegendentry{\hspace{0.5em}Approx. structured \panoc{} with improved line search}
        \legend{}
        \end{axis}
    \end{tikzpicture}
    \pgfplotstableread[col sep=comma]{data/warm10-2/inner_it.csv}\table
    \begin{tikzpicture}
        \begin{axis}[
            SmallBarPlot,
            xticklabels from table={\table}{label},
            ylabel=Number of iterations,
            xlabel=,
            width=0.45\linewidth,
            height=\myfigheight,
        ]
        \addplot [fill=pal1]
        table [x expr=\coordindex, y=s0] {\table};
        \addlegendentry{\hspace{0.5em}Original \panoc}
        \addplot [fill=pal2]
        table [x expr=\coordindex, y=s1] {\table};
        \addlegendentry{\hspace{0.5em}Original \panoc{} with improved line search}
        \addplot [fill=pal3]
        table [x expr=\coordindex, y=s2] {\table};
        \addlegendentry{\hspace{0.5em}Structured \panoc}
        \addplot [fill=pal4]
        table [x expr=\coordindex, y=s3] {\table};
        \addlegendentry{\hspace{0.5em}Structured \panoc{} with improved line search}
        \addplot [fill=pal5]
            table [x expr=\coordindex, y=s4] {\table};
        \addlegendentry{\hspace{0.5em}Approx. structured \panoc}
        \addplot [fill=pal6]
            table [x expr=\coordindex, y=s5] {\table};
        \addlegendentry{\hspace{0.5em}Approx. structured \panoc{} with improved line search}
        \legend{}
        \end{axis}
    \end{tikzpicture}
    \pgfplotstableread[col sep=comma]{data/warm10-2/runtime.csv}\table
    \begin{tikzpicture}
        \begin{axis}[
            SmallBarPlot,
            xticklabels from table={\table}{label},
            ylabel=Seconds,
            xlabel=,
            width=0.45\linewidth,
            height=\myfigheight,
        ]
        \addplot [fill=pal1]
            table [x expr=\coordindex, y=s0] {\table};
        \addlegendentry{\hspace{0.5em}Original \panoc}
        \addplot [fill=pal2]
            table [x expr=\coordindex, y=s1] {\table};
        \addlegendentry{\hspace{0.5em}Original \panoc{} with improved line search}
        \addplot [fill=pal3]
            table [x expr=\coordindex, y=s2] {\table};
        \addlegendentry{\hspace{0.5em}Structured \panoc}
        \addplot [fill=pal4]
            table [x expr=\coordindex, y=s3] {\table};
        \addlegendentry{\hspace{0.5em}Structured \panoc{} with improved line search}
        \addplot [fill=pal5]
            table [x expr=\coordindex, y=s4] {\table};
        \addlegendentry{\hspace{0.5em}Approx. structured \panoc}
        \addplot [fill=pal6]
            table [x expr=\coordindex, y=s5] {\table};
        \addlegendentry{\hspace{0.5em}Approx. structured \panoc{} with improved line search}
        \legend{}
        \end{axis}
    \end{tikzpicture}
    \ifMergefig
    \caption{Experimental results of the six solvers over 10 NMPC time steps 
    of the chain model, first without warm-starting (top) 
    and then warm-started using the shifted solution from the 
    previous time step (bottom).}
    \label{fig:cold-warm-start}
    \else
    \caption{Experimental results of the six solvers over 10 NMPC time steps 
    of the chain model, warm-started using the shifted solution from the 
    previous time step.}
    \label{fig:warm-start}
    \fi
\end{figure}

In most cases, the improved line search condition from \Cref{sec:impr-ls} reduces the 
number of function evaluations, the number of iterations, and the runtime. This 
is especially noticeable when the solver is cold-started.

When cold-started, the structured \panoc{} solver with the Hessian-vector product 
of equation \eqref{eq:61c} performs significantly better than standard \panoc{} 
in terms of the number of iterations. However, the run time is similar. The 
reason for this is the Hessian-vector term in \eqref{eq:61c}, which requires 
an additional gradient evaluation in each iteration when using finite differences
to compute this product. In the ideal case, standard \panoc{} requires 
one gradient evaluation per iteration, whereas structured \panoc{} with finite 
differences requires at least two, and this results in a higher run time per 
iteration.

The last two solvers approximate \eqref{eq:61c} by dropping the Hessian-vector 
term entirely. This eliminates one gradient evaluation per iteration, which 
can clearly be seen in the figures. The number of iterations is very close to
that of structured \panoc{} with the Hessian-vector term, so dropping this term 
results in a significant reduction of the run time.

When warm-starting, the differences in the number of iterations between the 
six solvers are much smaller, but structured \panoc{} without the Hessian-vector 
term still outperforms standard \panoc. When including the Hessian-vector term, 
the number of gradient evaluations increases by almost a factor of two, as 
described above, resulting in a significantly worse run time.

\subsection{Comparison to \ipopt{}}

In the following two experiments, the \alpaqa{} implementation of
structured \panoc{} with the improved line search condition without 
Hessian-vector products
(the best \panoc{}-based solver from the previous 
subsection) 
is compared to the popular interior point solver \ipopt{} \cite{ipopt}.

\begin{figure}[ht!]
    \centering
    \usepgfplotslibrary{colorbrewer}
    \pgfplotsset{every axis plot/.append style={mark size=1pt}}
    \tikzset{lbl/.style={anchor=west}}
    \begin{minipage}{0.48\columnwidth}
        \begin{tikzpicture}
            \begin{semilogyaxis}[
                font=\small,
                width=\linewidth,
                height=5cm,
                cycle list/Set1,
                cycle multiindex* list={
                    mark list*\nextlist
                    Set1\nextlist
                },
                title={Cold start},
                enlargelimits=false,
                ymin=0.1,
                ymax=12,
                xmin=0,
                ylabel=Solver run time {[s]},
                xlabel=MPC time step,
            ]
                \addplot+ table [x=x, y=Approx. structured PANOC with improved line search, col sep=comma] {data/ipopt-cold/ipopt-times.csv};
                \addplot+ table [x=x, y=Ipopt, col sep=comma] {data/ipopt-cold/ipopt-times.csv};
                \addplot+ table [x=x, y=Ipopt limited memory, col sep=comma] {data/ipopt-cold/ipopt-times.csv};
                \addplot+ table [x=x, y=Ipopt limited memory JIT, col sep=comma] {data/ipopt-cold/ipopt-times.csv};
            \end{semilogyaxis}
        \end{tikzpicture}
    \end{minipage}
    \begin{minipage}{0.48\columnwidth}
        \begin{tikzpicture}
            \begin{semilogyaxis}[
                font=\small,
                width=\linewidth,
                height=5cm,
                cycle list/Set1,
                cycle multiindex* list={
                    mark list*\nextlist
                    Set1\nextlist
                },
                title={Warm start},
                enlargelimits=false,
                ymin=0.1,
                ymax=12,
                xmin=0,
                yticklabels={},
                clip=false,
                xlabel=MPC time step,
            ]
                \addplot+ table [x=x, y=Approx. structured PANOC with improved line search, col sep=comma] {data/ipopt-warm/ipopt-times.csv} node[lbl]{\small \alpaqa};
                \addplot+ table [x=x, y=Ipopt, col sep=comma] {data/ipopt-warm/ipopt-times.csv} node[lbl]{\small \ipopt};
                \addplot+ table [x=x, y=Ipopt limited memory, col sep=comma] {data/ipopt-warm/ipopt-times.csv} node[lbl,anchor=south west]{\ipopt{} (\small \textsc{lm})};
                \addplot+ table [x=x, y=Ipopt limited memory JIT, col sep=comma] {data/ipopt-warm/ipopt-times.csv} node[lbl, anchor=north west, align=left]{\small \ipopt\\ (\textsc{lm+jit})};
            \end{semilogyaxis}
        \end{tikzpicture}
    \end{minipage}
    \caption{Comparison of the run time of \alpaqa{} (approximate structured \panoc{} with improved line search)
    and \ipopt{} (exact Hessian, limited memory Hessian, and limited memory with JIT compilation) 
    for 30 optimal control problems (with and without warm starting).}
    \label{fig:cold-warm-start-ipopt}
\end{figure}
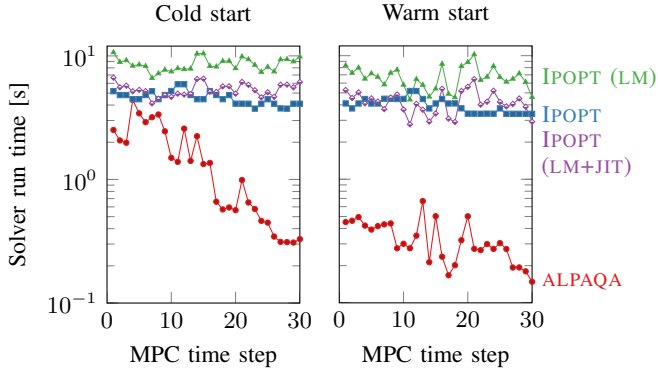

The tolerances used for \alpaqa{} are \(\varepsilon = \delta = 10^{-3}\).
For \ipopt{}, the \texttt{tol} and \texttt{constr\_viol\_tol}
options are both set to \(10^{-3}\). Three different configurations of 
\ipopt{} are used: (1) the default configuration with the exact Hessian
without just-in-time compilation\footnote{The \texttt{nlpsol jit} option 
was set to \texttt{False} in \casadi{}.}, 
(2) limited-memory Hessian 
approximation without just-in-time compilation, and 
(3) limited-memory Hessian approximation with just-in-time compilation%
\footnote{\ipopt{} with exact Hessian and just-in-time compilation is not included because it took at least 16 hours and 64 GiB of RAM to compile.}.
The four solvers are applied to the same hanging chain MPC problem as before
and simulated for 30 time steps. In a first experiment, the solvers are not 
warm-started (initial guesses set to zero), and in a second experiment, all 
solvers are warm-started using the shifted solution and multipliers from the previous time step.
The results are shown in 
\ifMergefig
\Cref{fig:cold-warm-start-ipopt}.
\else
\Cref{fig:cold-start-ipopt,fig:warm-start-ipopt}.
\fi

It is clear that the \alpaqa{} solver benefits from the warm start, being 
up to an order of magnitude faster compared to a cold start. 
\ipopt{} also 
converges slightly more quickly when warm-started, but the difference is not as 
substantial. Another observation is that the run time for \ipopt{} is relatively
consistent across time steps, whereas the performance of \alpaqa{} depends
on the constraints: in the first iterations of the simulation, many 
state constraints are active, and more ALM iterations are required to satisfy 
them, at the end of the simulation, this is no longer the case, and \alpaqa{}
converges more quickly.

For smaller tolerances, \ipopt{} with the exact Hessian 
does have an advantage over the first-order \alpaqa{} methods, but very precise 
solutions are usually not required in real-time MPC applications.

\subsection{Performance on the \cutest{} benchmarks}

\begin{figure}[ht!]
\centerline{
\includegraphics[width=0.495\linewidth]{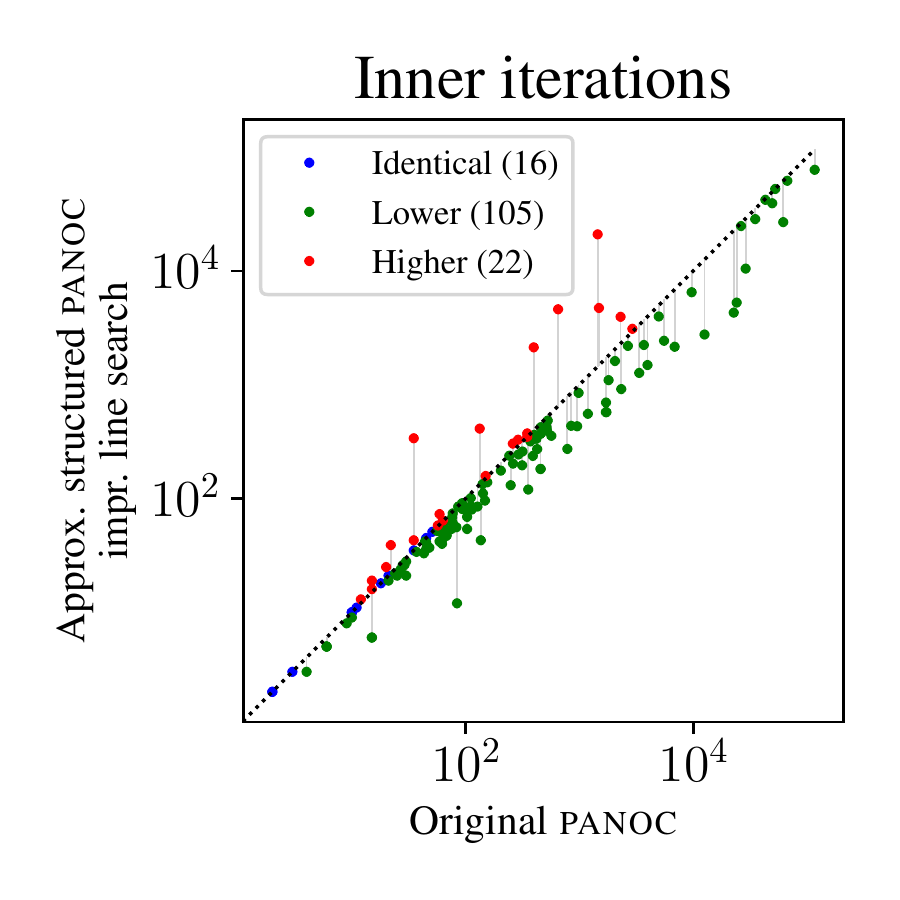}
\includegraphics[width=0.495\linewidth]{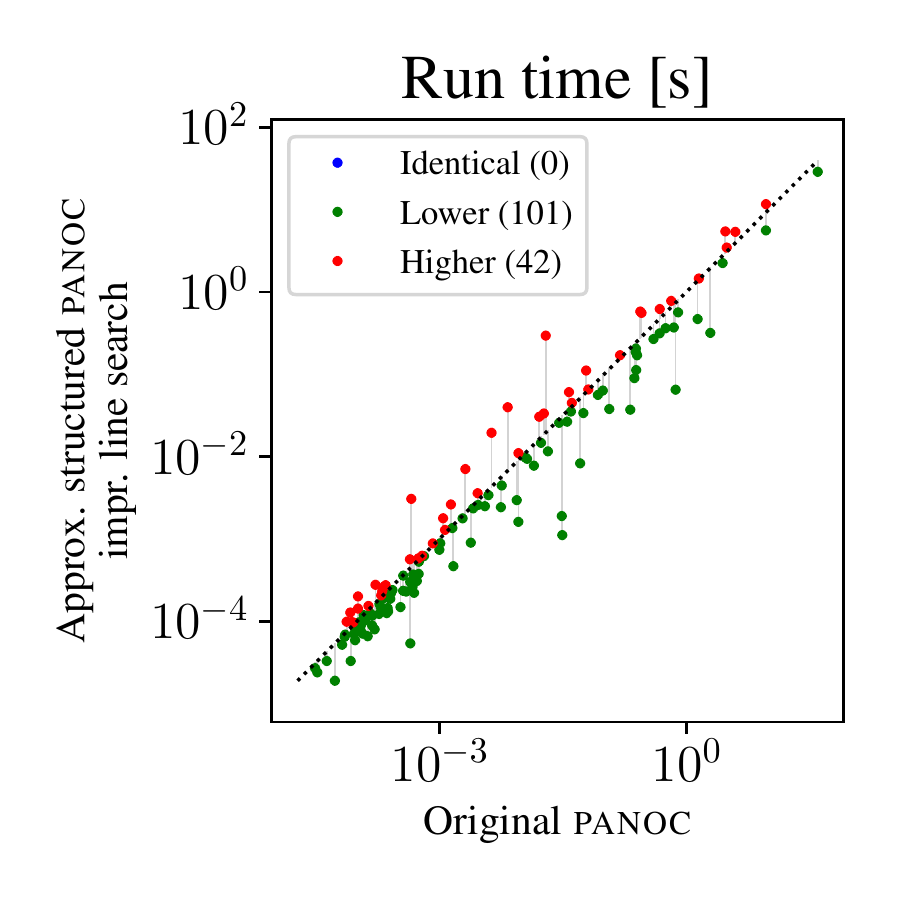}
}
\caption{Comparison between the number of iterations and run time for the 
original \panoc{} solver and the structured \panoc{} solver with the 
improved line search, applied to \cutest{} benchmarks.}
\label{fig:cutest}
\end{figure}

When applied to a collection of 219 problems (excluding QPs) from the 
\cutest{} benchmarks, the original \panoc{} solver solves 
148 problems, \panoc{} with the improved line search condition solves
153, and structured \panoc{} with the improved line search condition 
manages to solve 158. 
The latter is not only more robust, 
it also solves the problems more quickly than the original \panoc{}
solver, as shown in \Cref{fig:cutest}.

\section{Conclusion and further research}
We presented the \alpaqa{} library for nonconvex constrained 
optimization, and applied it to several benchmark problems. 
An attractive property of the augmented Lagrangian method and the \panoc{} 
algorithm is that it can be warm-started, making it competitive with 
state-of-the-art solvers such as \ipopt{} in MPC applications. The use of 
first-order methods opens the door to large-scale problems and embedded 
environments.

The \alpaqa{} library also implements two improvements to the original \panoc{} 
algorithm: a way to compute the quasi-Newton steps while exploiting the 
structure of the box-constrained inner problems, and a stricter line 
search condition that can reject low-quality quasi-Newton steps. 
These modifications were shown to improve both the performance 
and the robustness of \panoc{}.

Possible further work includes adding support for the quadratic penalty method
for handling non-smooth constraints, 
implementing the second-order solution strategies covered in 
\Cref{subsec:aughess-newton,subsec:lagrhess-newton}, and 
exploitation of the specific structure that arises from MPC problems.

\renewcommand*{\bibfont}{\footnotesize}
\printbibliography

@inproceedings{panoc,
  author        = {{Stella}, L. and {Themelis}, A. and {Sopasakis}, P. and {Patrinos}, P.},
  booktitle     = {2017 IEEE 56th Annual Conference on Decision and Control (CDC)},
  title         = {{A simple and efficient algorithm for nonlinear model predictive control}},
  year          = {2017},
  volume        = {},
  number        = {},
  pages         = {1939--1944},
  archiveprefix = {arXiv},
  eprint        = {1709.06487}
}

@book{nocedal_numopt,
  author    = {Nocedal, Jorge and Wright, Stephen J.},
  title     = {Numerical Optimization},
  year      = {2006},
  publisher = {Springer},
  address   = {New York}
}

@article{qpalm,
  title         = {{QPALM: A Proximal Augmented Lagrangian Method for Nonconvex Quadratic Programs}},
  author        = {Hermans, Ben and Themelis, Andreas and Patrinos, Panagiotis},
  year          = {2020},
  archiveprefix = {arXiv},
  eprint        = {2010.02653}
}

@article{Andersson2019,
  author    = {Joel A. E. Andersson and Joris Gillis and Greg Horn
               and James B Rawlings and Moritz Diehl},
  title     = {{CasADi} -- {A} software framework for nonlinear optimization
               and optimal control},
  journal   = {Mathematical Programming Computation},
  volume    = {11},
  number    = {1},
  pages     = {1--36},
  year      = {2019},
  publisher = {Springer}
}

@article{nocedal-lbfgs,
  title            = {{On the limited memory BFGS method for large scale optimization}},
  author           = {Liu, {Dong C.} and Jorge Nocedal},
  disable_note     = {Copyright: Copyright 2007 Elsevier B.V., All rights reserved.},
  year             = {1989},
  month            = aug,
  disable_doi      = {10.1007/BF01589116},
  disable_language = {English (US)},
  volume           = {45},
  pages            = {503--528},
  journal          = {Mathematical Programming},
  disable_issn     = {0025-5610},
  publisher        = {Springer-Verlag GmbH and Co. KG},
  number           = {1-3}
}

@article{ThemelisAndreas2018Feft,
  disable_issn  = {1052-6234},
  journal       = {Siam Journal On Optimization},
  pages         = {2274--2303},
  volume        = {28},
  publisher     = {Society for Industrial and Applied Mathematics},
  number        = {3},
  year          = {2018},
  title         = {{Forward-backward envelope for the sum of two nonconvex functions: further properties and nonmonotone line-search algorithms}},
  author        = {Themelis, Andreas and Stella, Lorenzo and Patrinos, Panos},
  archiveprefix = {arXiv},
  eprint        = {1606.06256}
}

@article{dennismore,
  oissn     = {00255718, 10886842},
  ourl      = {http://www.jstor.org/stable/2005926},
  author    = {J. E. Dennis and Jorge J. Moré},
  journal   = {Mathematics of Computation},
  number    = {126},
  pages     = {549--560},
  publisher = {American Mathematical Society},
  title     = {{A Characterization of Superlinear Convergence and Its Application to Quasi-Newton Methods}},
  volume    = {28},
  year      = {1974}
}

@book{practicalaugmentedlagrangian,
  author    = {E. G. Birgin and J. M. Martínez},
  editor    = {Nicholas J. Higham},
  title     = {Practical Augmented Lagrangian Methods for Constrained Optimization},
  publisher = {SIAM},
  year      = {2014}
}

@inproceedings{PatrinosPanagiotis2013PNmf,
  disable_issn = {0191-2216},
  pages        = {2358--2363},
  publisher    = {IEEE},
  booktitle    = {52nd IEEE Conference on Decision and Control},
  disable_isbn = {1467357146},
  year         = {2013},
  title        = {{Proximal Newton methods for convex composite optimization}},
  author       = {Patrinos, Panagiotis and Bemporad, Alberto}
}

@book{RockafellarVariationalAnalysis,
  series       = {Grundlehren der mathematischen Wissenschaften 317},
  publisher    = {Springer},
  disable_isbn = {3540627723},
  year         = {2004},
  title        = {{Variational analysis}},
  address      = {Berlin},
  author       = {Rockafellar, R. Tyrrell and Wets, Roger J.-B.}
}

@article{ipopt,
  author       = {W\"achter, Andreas and Biegler, Lorenz T.},
  title        = {{On the implementation of an interior-point filter line-search algorithm for large-scale nonlinear programming}},
  journal      = {Mathematical Programming},
  volume       = {106},
  number       = {1},
  year         = {2006},
  pages        = {25--57},
  disabled_url = {https://doi.org/10.1007/s10107-004-0559-y}
}

@article{BensonHandeY2008Imfn,
  disabled_issn = {0926-6003},
  journal       = {Computational optimization and applications},
  pages         = {143--189},
  volume        = {40},
  publisher     = {Springer US},
  number        = {2},
  year          = {2008},
  title         = {Interior-point methods for nonconvex nonlinear programming: regularization and warmstarts},
  copyright     = {Springer Science+Business Media, LLC 2008},
  address       = {Boston},
  author        = {Benson, Hande Y and Shanno, David F}
}

@article{FangHaw-ren2008MCaa,
  disable_issn = {0025-5610},
  journal      = {Mathematical programming},
  pages        = {319--349},
  volume       = {115},
  publisher    = {Springer-Verlag},
  number       = {2},
  year         = {2008},
  title        = {{Modified Cholesky algorithms: a catalog with new approaches}},
  copyright    = {Springer-Verlag 2007},
  address      = {Berlin/Heidelberg},
  author       = {Fang, Haw-Ren and O’Leary, Dianne P}
}

@article{DuffIain2004Mcft,
  disable_issn = {0098-3500},
  journal      = {ACM transactions on mathematical software},
  pages        = {118--144},
  volume       = {30},
  publisher    = {ACM},
  number       = {2},
  year         = {2004},
  title        = {{MA57---a code for the solution of sparse symmetric definite and indefinite systems}},
  copyright    = {Copyright 2012 Elsevier B.V., All rights reserved.},
  address      = {NEW YORK},
  author       = {Duff, Iain}
}

@article{open,
  journal       = {IFAC-PapersOnLine},
  volume        = {53},
  number        = {2},
  pages         = {6548-6554},
  year          = {2020},
  note          = {21st IFAC World Congress},
  disabled_issn = {2405-8963},
  disabled_doi  = {https://doi.org/10.1016/j.ifacol.2020.12.071},
  disabled_url  = {https://www.sciencedirect.com/science/article/pii/S240589632030327X},
  author        = {Pantelis Sopasakis and Emil Fresk and Panagiotis Patrinos},
  title         = {{OpEn: Code Generation for Embedded Nonconvex Optimization}},
  archiveprefix = {arXiv},
  eprint        = {2003.00292}
}

@inproceedings{chainmodel,
  author      = {Wirsching, Leonard and Bock, Hans and Diehl, Moritz},
  year        = {2006},
  month       = {11},
  pages       = {591 - 596},
  title       = {{Fast NMPC of a chain of masses connected by springs}},
  journal     = {Proceedings of the IEEE International Conference on Control Applications},
  disable_doi = {10.1109/CACSD-CCA-ISIC.2006.4776712}
}

@inproceedings{8550253,
  author        = {Sathya, Ajay and Sopasakis, Pantelis and Van Parys, Ruben and Themelis, Andreas and Pipeleers, Goele and Patrinos, Panagiotis},
  booktitle     = {2018 European Control Conference (ECC)},
  title         = {{Embedded nonlinear model predictive control for obstacle avoidance using PANOC}},
  year          = {2018},
  volume        = {},
  number        = {},
  pages         = {1523-1528},
  disabled_doi  = {10.23919/ECC.2018.8550253},
  archiveprefix = {arXiv},
  eprint        = {1904.10546}
}

@inproceedings{small_AerialNavigationObstructed_2019,
  title         = {{Aerial Navigation in Obstructed Environments with Embedded Nonlinear Model Predictive Control}},
  booktitle     = {2019 18th {{European Control Conference}} ({{ECC}})},
  author        = {Small, Elias and Sopasakis, Pantelis and Fresk, Emil and Patrinos, Panagiotis and Nikolakopoulos, George},
  year          = {2019},
  month         = jun,
  pages         = {3556--3563},
  disabled_doi  = {10.23919/ECC.2019.8796236},
  archiveprefix = {arXiv},
  eprint        = {1812.04755}
}

\end{document}